\documentclass[runningheads]{llncs}
\usepackage[font=small]{subfig}
\usepackage[compatible]{algpseudocode}
\usepackage{algorithm}
\usepackage{graphicx}
\usepackage{epstopdf}
\usepackage{multirow}
\begin{document}
\title{Metaheuristics for Min-Power Bounded-Hops Symmetric Connectivity Problem\thanks{ The research is supported by the Russian Science Foundation (project 18-71-00084).}}
\titlerunning{Metaheuristics for MPBHSCP}
%
%

\author{Roman Plotnikov\inst{1}\orcidID{0000-0003-2038-5609} \and
Adil Erzin\inst{1,2}\orcidID{0000-0002-2183-523X}}
\authorrunning{R. Plotnikov et al.}
%

\institute{Sobolev Institute of Mathematics, Novosibirsk, Russia \and
Novosibirsk State University, Novosibirsk, Russia}

\maketitle              

The article is submitted to the special session: Synergy of machine learning, combinatorial optimization, and computational geometry: problems, algorithms, bounds, and applications.

\begin{abstract}
We consider a Min-Power Bounded-Hops Symmetric Connectivity problem that consists of the construction of communication spanning tree on a given graph, where the total energy consumption spent for the data transmission is minimized and the maximum number of edges between two nodes is bounded by some predefined constant. We focus on the planar Euclidian case of this problem where the nodes are placed at the random uniformly spread points on a square and the power cost necessary for the communication between two network elements is proportional to the squared distance between them. Since this is an NP-hard problem, we propose different heuristics based on the following metaheuristics: genetic local search, variable neighborhood search, and ant colony optimization. We perform a posteriori comparative analysis of the proposed algorithms and present the obtained results in this paper.

\keywords{Energy efficiency  \and Approximation algorithms \and Symmetric connectivity \and Bounded hops \and Genetic local search \and Variable neighborhood search \and Ant colony optimization.}

\end{abstract}

\section{Introduction}
Due to the prevalence of wireless sensor networks (WSNs) in human life, the different optimization problems aimed to increase their efficiency remain actual. Since usually WSN consists of elements with the non-renewable power supply with restricted capacity, one of the most important issues related to the design of WSN is prolongation its lifetime by minimizing energy consumption of its elements per time unit. A significant part of sensor energy is spent on the communication with other network elements. Therefore, the modern sensors often have an ability to adjust their transmission ranges changing the transmitter power. Herewith, usually, the energy consumption of a network's element is assumed to be proportional to $d^s$, where $s\ge 2$ and $d$ is the transmission range \cite{R96}.

The problem of determining the optimal power assignment in WSN is well-studied. The most general Range Assignment Problem, where the goal is to find a strongly connected subgraph in a given directed graph, has been considered in \cite{CPS99,KKKP00}. Its subproblem, Min-Power Symmetric Connectivity Problem (MPSCP), was first studied in \cite{CMZ02}. The authors proved that Minimum Spanning Tree (MST) is a 2-approximation solution to this problem. Also, they proposed a polynomial-time approximation scheme with a performance ratio of $1 + \ln{2} + \varepsilon \approx 1.69$ and a 15/8-approximation polynomial algorithm. In \cite{CNSCL03} a greedy heuristic, later called  Incremental Power: Prim (IPP), was proposed. The IPP is similar to the Prim's algorithm for MST constructing. A Kruskal-like heuristic, later called Incremental Power: Kruskal, was studied in \cite{CN02}. Both of these so-called incremental power heuristics have been proposed for the Minimum Power Asymmetric Broadcast Problem, but they are suitable for MPSCP too. It is proved in \cite{PS06} that they both have an approximation ratio 2, and it was shown in the same paper that in practice they yield significantly more accurate solution than MST. Also, in a series of papers different heuristic algorithms have been proposed for MPSCP and the experimental studies have been done: local search procedures \cite{PS06,ACMPTZ06,EPS13}, methods based on iterative local search \cite{WM09}, hybrid genetic algorithm that uses a variable neighborhood descent as mutation \cite{EP15}, variable neighborhood search \cite{EMP16_COR}, and variable neighborhood decomposition search \cite{PEM18_VNDS}.

Another important property of WSN's efficiency is a message transmission delay, i.e., the minimum time necessary for transmitting a message from one sensor to another via the intermediate transit nodes. As a rule, the delay is proportional to the maximum number of hops (edges) between two nodes of a network. In the general case, when the network is represented as a directed arc-weighted graph, and the goal is to find a strongly connected subgraph with minimum total power consumptions and bounded path length, the problem is called a Min-Power Bounded-Hops Strong Connectivity Problem. In \cite{Clementi00_1} the approximation algorithms with guaranteed estimates have been proposed for the Euclidean case of this problem. The bi-criteria approximation algorithm for the general case (not necessarily Euclidian) with guaranteed upper bounds has been proposed in \cite{Calinescu06}. The authors of \cite{Carmi15} propose an improved constant factor approximation for the planar Euclidian case of the problem.

In this paper, we consider the symmetric case of Min-Power Bounded-Hops Strong Connectivity Problem, when the network is represented as an undirected edge-weighted graph. Such a problem is known as Min-Power Bounded-Hops Symmetric Connectivity Problem (MPBHSCP) \cite{Calinescu06}. We also assume that the sensors are positioned on Euclidian plane. The energy consumption for the data transmission is assumed to be proportional to the area of a circle with center in sensor position and radius equal to its transmission range $d$, and, therefore, $s = 2$. This problem is still NP-hard in planar Euclidian case \cite{Clementi00}, and, therefore, the approximation heuristic algorithms that allow obtaining the near-optimal solution in a short time, are required for it.

A set of polynomial algorithms that construct the approximate solutions to MPBHSCP were proposed in \cite{PloEr19}. In this pasper, we suggest three metaheuristic approaches that aimed to improve the solutions obtained by the known constructive heuristics. Namely, we use a variable neighborhood search, a genetic local search, and an ant colony optimization that make use of different variants of local search procedure. This research was inspired by the papers where the different metaheuristics are successfully applied for the approximate solution of Bounded-Diameter Minimum Spanning Tree (BDMST) (e.g., see \cite{Gruber06}), and for MPSCP (e.g., see \cite{EMP16_COR}). We conducted an extensive numerical experiment to compare our algorithms. We present the results of the experiment in this paper. Note that, to the best of our knowledge, the metaheuristics of such kind previously were never applied to MPBHSCP.

The rest of the paper is organized as follows. In Section \ref{sPF} the problem is formulated, in Section \ref{sH} descriptions of the proposed algorithms are given, Section \ref{sS} contains results and analysis of the experimental study, and Section \ref{sC} concludes the paper.

\section{Problem formulation}\label{sPF}

Mathematically, MPBHSCP can be formulated as follows. Given a connected edge-weighted undirected graph $G = (V,E)$ and an integer value $D \geq 1$, find such spanning
tree $T^*$ in $G$, which is the solution to the following problem:

\begin{equation}\label{e1}
    W(T)=\sum\limits_{i\in V}\max\limits_{j\in V_i(T)}
    c_{ij}\rightarrow\min\limits_T,
\end{equation}

\begin{equation}\label{e2}
    dist_T(u, v) \leq D \; \forall u,v \in V,
\end{equation}

where $V_i(T)$ is the set of vertices adjacent to the vertex $i$ in the
tree $T$, $c_{ij} \geq 0$ is the weight of the edge $(i,
j) \in E$, and $dist_T(u, v)$ is the number of edges in a path between the vertices $u \in V$ and $v \in V$ in $T$.

Obviously, in general case, MPBHSCP may even not have any feasible solution. In this paper, we consider a planar Euclidian case, where an edge weight equals the squared distance between the corresponding points and $G$ is a complete graph. Therefore, a solution always exists.

Although any feasible solution of (\ref{e1})--(\ref{e2}) is an undirected spanning tree with bounded diameter, we always can choose a center of this tree, i.e., a vertex (or two vertices if $D$ is odd), such that a path from it to any other vertex in a tree contains not more than $\lfloor D/2 \rfloor$ edges. Therefore, it is convenient to consider a solution as a directed tree (or arborescence) rooted in one of its centers. Further, we assume that the centers and the root are predefined for each considered feasible spanning tree, and, therefore, we will handle with the following notations that are suitable for directed trees: $v_0$ --- a root of a tree $T = (V, E_T)$; $P_T(v)$ --- a parent vertex of $v \in V \setminus \{v_0\}$ in $T$; $L_T(v)$ --- a level (i.e., the number of edges in a path from $v$ to the center) of $v \in V$ in $T$.

\section{Heuristic algorithms}\label{sH}
In this section we will describe the heuristic algorithms for approximation solution of MPBHSCP. Our methods are based on the following metaheuristics: variable neighborhood search (VNS), genetic local search (GLS), and ant colony optimization (ACO). All of our methods start with some initial feasible solution -- a spanning tree with bounded diameter (or with a set of feasible solutions, as in GLS). We assume that at least one such solution was already constructed by some heuristic, and the goal of our algorithms is to improve it in a best possible way.

Besides the obvious differences that are specific for the particular metaheuristics, our algorithms have the common parts. Namely, they use the same variants of local search and random movement procedures. Therefore, we will describe these procedures at first.

\subsection{Local search}\label{ssLS}
We suggest three types of neighborhood structure that are used in the local search procedures. The first neighborhood movement is called \emph{LevelChange}. It consists of changing a parent node for a vertex in a such way that the level of a vertex changes and the diameter is feasible (at most $D$). The second procedure, \emph{SameLevelParentChange}, consists of changing a parent for a vertex preserving its level. And the last one is \emph{CenterChange}, which consists of changing of a center vertex of a tree. Note that none of these three variants of local movement can be replaced by a sequence of others.

The first two local movements, \emph{LevelChange} and \emph{SameLevelParentChange} are quite simple. In both cases, at first, one edge $e = (v, P_T(v))$ is removed from $T$, and then another vertex $v_1$, that is not a descendant of $v$ in $T$, is chosen as a new parent of $v$. Herewith, some special conditions should be met: in the case of \emph{LevelChange}, $L_T(v_1)$ should not be equal to $L_T(v) - 1$ and the diameter restriction should not be violated; in the case of \emph{SameLevelParentChange}, the equality $L_T(v_1) = L_T(v) - 1$ should hold.

In the \emph{CenterChange} movement, at first, one center $c \in V$ is chosen (it may be either a root or another center in a case of odd diameter), and then, some other non-center vertex $v \in V$ is chosen as a new center. In order to make $v$ a new center instead of $c$ the following steps are performed: (a) the children of $c$ change their parent from $c$ to $v$; (b) $v$ is detached from its parent $v_p = P_T(v)$; (c) if $c$ is a root then $v$ becomes a root, otherwise it becomes a second center, and the root $v_0$ becomes a parent of $v$; (d) if $c \neq v_p$, then $c$ becomes a child of $v_p$, otherwise, it becomes a child of $v$.

Our algorithms use these three variants of neighborhood movement as parts of one local search method based on variable neighborhood descent metaheuristic (VND). The idea of VND is to perform local search within more than one neighborhood structure. This approach was proposed in \cite{HansenVNS}, and the pseudo-code is given in Algorithm \ref{alg:vnd}. In result, VND returns a local optimum for all considered neighborhood structures.

\begin{algorithm}[!bp]
\begin{algorithmic}[1]
\STATE Select an initial solution $T$;
\STATE $k \leftarrow 0$;
\STATE Set the set of the local searches $(LS_{l})_{l=1, 2, 3} $ $ \leftarrow $ $\{$\emph{LevelChange}, \emph{SameLevelParentChange}, \emph{CenterChange}$\}$;
\STATE $improved \leftarrow $ \textbf{true};
\WHILE {$improved$}
\STATE $improved \leftarrow $ \textbf{false}, $l \leftarrow 1$;
\WHILE {$l \leq 3$}
\STATE $T^{\prime} \leftarrow LS_l(T)$;
\IF {$T^{\prime}$ is better than $T$}
\STATE $T \leftarrow T^{\prime}$, $l \leftarrow 1$, $improved \leftarrow $ \textbf{true};
\ELSE
\STATE $l \leftarrow l+1$;
\ENDIF

\ENDWHILE
\ENDWHILE

\end{algorithmic}
\caption{Variable neighborhood descent} \label{alg:vnd}
\end{algorithm}

\subsection{Random movement}
Besides the local search procedures, some of our metaheuristics (to be precise, GLS and VNS) involve an operator of randomized modification of a tree. For such random movement we suggest a procedure \emph{RandomBranchReattaching}. In this procedure, some edge $(v, P_T(v))$ is chosen at random and removed from $T$. Then, $v$ is connected with a non-descendant vertex $u$, which is chosen at random as well, if this operation keeps the feasibility of a tree. This process is repeated $k$ times, where $k$ is an external integer parameter provided by an upper-level metaheuristic.

\subsection{Variable neighborhood search}
Variable neighborhood search (VNS) is a metaheuristic developed by Hansen and Mladenovic \cite{HansenVNS}, and it consists of two phases: randomized phase, or so-called shaking, when the current solution  is changed in random or in half-random way, and deterministic phase, where VND is applied to the shaken solution. In our implementation, \emph{RandomBranchReattaching} is used for the shaking phase, and the neighborhood movements \emph{LevelChange}, \emph{SameLevelParentChange}, and \emph{CenterChange} are used in the local search phase. The pseudo-code is presented in Algorithm \ref{alg:vns}. The great advantage of this metaheuristic, comparing to others, is that it requires tuning of the only parameter $k_{\max}$. The algorithm starts with some feasible solution. As the first approximation for MPBHSCP we use the best of the trees obtained by the heuristic algorithms that are proposed in \cite{PloEr19}: MPCBTC, MPRTC, MPCBLSoC, MPCBRC,MPQBH, and MPIR.

\begin{algorithm}[h]
\begin{algorithmic}[1]
\STATE Select an initial solution $T$;
\STATE $k \leftarrow 0$;
\WHILE {the stopping criteria is not met}
\WHILE {$k \leq k_{\max}$}
\STATE Perform shaking: $T^{\prime} \leftarrow \emph{RandomBranchReattaching}(T, k)$;
\STATE Apply the local search procedures to the shaken solution: $T^{\prime\prime} \leftarrow VND(T^{\prime})$;
\IF {$T^{\prime\prime}$ is better than $T$}
\STATE $T \leftarrow T^{\prime\prime}$; $k \leftarrow 1$;
\ELSE
\STATE $k \leftarrow k+1$;
\ENDIF
\ENDWHILE
\ENDWHILE
\end{algorithmic}
\caption{Variable neighborhood search} \label{alg:vns}
\end{algorithm}

\subsection{Genetic local search}
Another approach suitable for the problem (\ref{e1})--(ref{e2}) is genetic local search algorithm. This metaheuristic deals with \emph{population} --- a set of feasible solutions. Before the algorithm starts, its first population should be generated. For the first population, we used all spanning trees constructed by 6 algorithms from \cite{PloEr19}. Note that, since MPRTC is randomized, it may construct a set of different feasible solutions instead of only one solution, that yield other 5 deterministic heuristics. This allows us to generate the first population with the size which does not exceed some predefined value. Each iteration of the algorithm consists of applying the following operators to the current population: (a) calculation of \emph{fitness}, that expresses the quality of the solution; (b) \emph{selection}, that chooses a subset of solutions from the population according to their fitness; (d) \emph{crossover}, that creates a new solution (an offspring) from the selected pair of solutions; (e) \emph{mutation}, that randomly modifies the offspring; (f) \emph{local search}, that improves the offspring; (g) \emph{join}, that selects the population of the next generation from the current population and the set of the offsprings. A brief description of the main steps of this algorithm is presented in Algorithm \ref{alg:gls}

\begin{algorithm}[!hbtp]
\begin{algorithmic}[1]
\STATE Generation of the first population;
\STATE Fitness calculation of the population;
\WHILE {stop condition is not met}
\STATE Selection;
\STATE Crossover;
\STATE Mutation;
\STATE Local search by VND;
\STATE Fitness calculation of the offspring;
\STATE Join;
\STATE $T \leftarrow$ the best tree among the current population;
\ENDWHILE
\end{algorithmic}
\caption{Genetic local search} \label{alg:gls}
\end{algorithm}

In our implementation of genetic local search, we take the value of $1 / W(T)$ as a fitness of $T$. This corresponds to the rule that fitness has to be a positive value which is higher when the value of the objective function is closer to optimum. Within the selection procedure, a set of prospective parents of the next offspring is filled with solutions from the current population in the following way. Sequentially, two trees are taken from the current population in proportion to their fitness probability: the first tree of each pair is chosen randomly from the entire population, and the second tree is chosen from the remaining part of the population. Each pair should contain different trees, but the same tree may be included in many pairs.

For the crossover operator, a solution is represented as an array of integer values, that correspond to the vertex levels in a tree. In other words, we assume that the vertices are numbered, and for each number  $i = 1, ..., n$, the value of $i$-th element in array is assigned to the level of $i$-th vertex in a tree. Given two integer arrays (let's call them the \emph{parent} arrays), a new (\emph{child}) array, that will correspond to an offspring, is generated in the following way. First of all, an offspring has to have a center. For that reason, one parent array is taken at random (with probability of 0.5), and then, the child array derives the elements assigned to 0 from this parent array. Each parent array has one or two such elements, depending on parity of $D$, and the child array should have the same number of elements assigned to 0. The values of all other elements in the child array derive the values at the same places from the parents, and each time the parent is chosen with probability of 0.5. Note that if the element that assigned to 0 is chosen to be derived by a child, then the corresponding element of a child is assigned to 1. This is done because the corresponding vertex cannot be a center of the offspring, since its center is already established.

The decoding of a tree from the integer array is performed in the following way. Let $A$ be the array of integers that should be decoded to a tree $T$. At first, a such vertex $v_0$, that $A(v_0) = 0$, is assigned to the root of a tree, and, if another vertex $v_1$ with the same property exists, then it is assigned to the second center of a tree and $v_0$ is assigned to the parent of $v_1$. After that, for each other $i$-th element of an array its predecessor in $T$ $j$ is chosen in such way, that $A(j) < A(i)$ and the edge that connects $i$-th and $j$-th vertices, brings the the minimum contribution to the value of the objective function.

The mutation procedure takes as an argument (an integer parameter) $k$ --- the maximum difference (number of different arcs in the initial tree and in the modified one). This parameter is taken randomly from the interval $[1, n / 3]$, with probability proportional to its inverse value (i.e., smaller modifications are more possible). To perform a random movement for the mutation, we used the procedure \emph{RandomBranchReattaching}. The mutation procedure is applied with probability $PM$ (a parameter of the algorithm) to each offspring.

Additionally, our algorithm applies local search to improve the offsprings after the crossover operator. To do this, we used VND algorithm that performs local search within three neighborhood structures defined above: \emph{LevelChange}, \emph{SameLevelParentChange}, and \emph{CenterChange}. This solution improvement procedure is applied with a predefined probability, as well as a randomized mutation.

At the \emph{join} procedure a subset of solutions from the current population and the current offspring, which have the largest fitness values, are chosen to fill the population of the next generation.

Our version of GLS for MPBHSCP requires the following parameters:

\begin{itemize}
  \item $PopSize$ --- the size of population;
  \item	$OffspSize$ --- the size of offspring;
  \item	$PM$ --- the probability of mutation;
  \item	$PLS$ --- the probability of local search.
\end{itemize}

\subsection{Ant colony optimization}
As the third heuristic algorithm for the approximate solution of MPBHSCP, we propose an algorithm based on the \emph{ant colony optimization metaheuristic} (ACO). A \emph{path} of an ant corresponds to the solution to the problem. The path usually consists of the elements, each of which is chosen randomly with probability depending on \emph{pheromone value} that stores information about the frequency of usage of a particular part of a path in the best-found solution. We designed our algorithm in a similar manner as it was done by Gruber et al. in \cite{Gruber06}. To represent a feasible solution of MPBHSCP as a path we used the same vertex-level encoding that was used in the crossover operator of GLS, i.e., an array of $n$ integers not greater than $\lfloor D/2 \rfloor$ corresponding to the vertex levels. As the pheromone values we used the matrix $(\tau_{il})$ of size $n \times \lfloor D/2 \rfloor$, that is initially filled with equal non-negative real numbers $1 / ( n \cdot W(T_0))$, where $T_0$ is the initial solution. Our variant of the ACO algorithm consists of three phases: (a) paths construction, (b) solutions improvement, and (c) pheromone matrix updating. The main steps of the algorithm are briefly described in Algorithm \ref{alg:aco}.

\begin{algorithm}[!hbtp]
\begin{algorithmic}[1]
\STATE Generation of pheromone matrix;
\WHILE {stop condition is not met}
\STATE Construction of ant paths according to the pheromone matrix;
\STATE Improvement of the solutions that are derived from the paths;
\STATE Update of the pheromone matrix;
\ENDWHILE
\end{algorithmic}
\caption{Ant colony optimization} \label{alg:aco}
\end{algorithm}

In the paths construction phase, at first, the center of a corresponding tree should be defined. For that reason, we assign one or two (depending on parity of $D$) elements of ant path to 0. The vertices (or indices of ant path) that will assigned to the centers (with level 0) are chosen randomly with the probability $P_{i,0} = \tau_{i, 0} / \sum_{j = 1}^{n}{\tau_{j, 0}}$. After that, for each vertex $i = 1, ..., n$, that has not been assigned to the center, its level is assigned randomly with the probability $P_{i, l} =  \tau_{i, l} / \sum_{l^{\prime} = 1}^{\lfloor D/2 \rfloor}{\tau_{j, l^{\prime}}}$, where $l = 1, ..., \lfloor D/2 \rfloor$.

After this construction, the paths are transformed into the spanning trees by the same decoding procedure that is used in the GLS. Each spanning tree is then improved by the VND procedure, that makes use of three neighborhood search types: \emph{LevelChange}, \emph{SameLevelParentChange}, and \emph{CenterChange}. After that, the best solution found so far $T_{best}$ is used for the updating of the pheromone matrix in the following way. For each $i = 1, ..., n$, and $l = 0, ..., \lfloor D/2 \rfloor$, $\tau_{i, l} = \tau_{i, l} + \rho / W(T_{best})$ if $l = L_{T_{best}}(i)$, and  $\tau_{i, l} = \tau_{i, l}(1 - \rho)$, otherwise. ACO requires two parameters: $ColSize$ --- the number of ants in colony, and $\rho$ --- the pheromone decay coefficient.

\section{Simulation} \label{sS}
We have implemented all the described algorithms in C++ programming language and launched them on the Intel Core i5-4460 3.2GHz processor with 8Gb RAM. In order to make our experiment results reproducible, we used as test instances the data sets that are given in Beasley's OR-Library for Euclidian Steiner Problem (http://people.brunel.ac.uk/~mastjjb/jeb/orlib). These test cases present the random uniformly distributed points in the unit square. We tested 3 variants of dimension: $n = $ 100, 250, and 500. We also took different values of $D$ for each dimension. Since all of our algorithms are partially probabilistic, we launched each algorithm 10 times on each instance, and calculated the average value of objective, the best value of objective, and the standard deviation. As a stop criteria the following condition was used: the best found solution is not changed during three iterations in a row.

We have performed the preliminary testing of each algorithm to determine such combination of its parameters that would provide in most cases the best result without consuming much time for the calculations. For VNS, $k_{\max}$ was chosen from the set $\{20, 30, 40\}$, and 30 appeared to be the best variant. For GLS, the pair $(PopSize, OffspSize) = (75, 40)$ appeared to be the best among the variants $\{ (25, 15)$, $(50, 20)$, $(75, 40)$, $(100, 50) \}$, and the pair $(PM, PLS)$ = $(0.5, 0.5)$ appeared to be the best among the variants $\{(0.25, 025)$, $(0.25, 05)$, $(0.25, 075)$, $(0.5, 05)$, $(0.75, 0.25)$, $(0.75, 075)\}$. As for ACO, we found out that $\rho = 0.2$ is the best choice among $\{ 0.005$, $0.01$, $0.05$, $0.1$, $0.2\}$ and $ColSize = 50$ is the best choice among $\{25$, $50$, $100\}$. We also tried to exclude one or more variants of local search in the VND subroutine of our algorithms, but, on average, this always deteriorated the results. Therefore, we decided to keep all the proposed variants of local search in each of our algorithms.

The results of the experiment are presented in Table \ref{tab:res}. The first three columns contain test instance properties: the tree diameter bound, $D$, the size of a problem, $n$, and the instance case number in the OR Library, nr. In the fourth column, the objective values on the best of constructive heuristics results, $T_{CH}$, are presented. Note that $T_{CH}$ was passed to each metaheuristic algorithm as the initial solution. We added this column here because it is important to see, how much the best solution found by constructive heuristics was improved by our metaheuristic based algorithms. In the other columns, the results of ACO, GLS, and VNS are presented: the objective values on the best found solutions, $W_{best}$, the average values of objective, $W_{av}$, the standard deviation of the set of objective values on the found solutions, $W_{sd}$, and the average running times. The best values among all algorithms are marked bold.

\begin{table}[!htp]
\begin{tabular}{|c|c|c|c|ccc|ccc|ccc|ccc|}
\hline

\multirow{2}{*}{$D$} &\multirow{2}{*}{$n$} &\multirow{2}{*}{nr} &\multirow{2}{*}{$W(T_{CH})$} &\multicolumn{3}{c|}{$W_{best}$} &\multicolumn{3}{c|}{$W_{av}$} &\multicolumn{3}{c|}{$W_{sd}$} &\multicolumn{3}{c|}{Time (in sec.)}\\

\cline{5-16}

& & &
&ACO &GLS &VNS
&ACO &GLS &VNS
&ACO &GLS &VNS
&ACO &GLS &VNS\\

\hline

\multirow{9}{*}{7}	&\multirow{3}{*}{50}
    &1  &1.89	&1.56	&\textbf{1.47}	&1.61	&1.63	&1.69	&\textbf{1.61}	&0.07	&0.18	&0	    &1.87	&0.69	&\textbf{0.36}\\
&	&2	&1.77	&\textbf{1.31}	&1.42	&1.34	&\textbf{1.39}	&1.55	&1.42	&0.06	&0.12	&0.04	&2.22	&\textbf{0.69}	&0.75\\
&	&3	&1.71	&\textbf{1.24}	&1.30	&1.36	&\textbf{1.33}	&1.44	&1.40	&0.05	&0.12	&0.02	&1.57	&\textbf{0.67}	&1.10\\

\cline{2-16}
&\multirow{3}{*}{100}

    &1  &2.07	&\textbf{1.60}	&1.74	&1.66	&1.96	&1.80	&\textbf{1.70}	&0.14	&0.05	&0.04	&5.93	&3.26	&\textbf{1.93}\\
&	&2	&2.00	&\textbf{1.48}	&1.55	&1.70	&1.73	&1.82	&\textbf{1.70}	&0.12	&0.13	&0.01	&8.16	&3.03	&\textbf{1.88}\\
&	&3	&2.35	&2.35	&1.99	&\textbf{1.90}	&2.35	&2.03	&\textbf{1.92}	&0	    &0.04	&0.04	&\textbf{2.60}	&3.15	&3.04\\

\cline{2-16}
&\multirow{3}{*}{250}
    &1  &3.13	&3.13	&2.76	&\textbf{2.60}	&3.13	&2.92	&\textbf{2.62}	&0	&0.09	&0.02	&28.82	&33.07	&\textbf{21.16}\\
&	&2	&3.30	&3.30	&2.80	&\textbf{2.88}	&3.30	&2.98	&\textbf{2.89}	&0	&0.17	&0.01	&28.73	&32.36	&\textbf{12.50}\\
&	&3	&3.11	&3.11	&2.53	&\textbf{2.43}	&3.11	&2.76	&\textbf{2.49}	&0	&0.15	&0.05	&28.57	&32.64	&\textbf{30.41}\\

\hline

\multirow{10}{*}{10}	&\multirow{3}{*}{50}
    &1  &1.68	&\textbf{1.14}	&1.23	&1.22	&1.29	&1.25	&\textbf{1.24}	&0.07	&0.03	&0.03	&1.78	&\textbf{0.76}	&1.05\\
&	&2	&1.18	&\textbf{1.01}	&1.13	&1.06	&1.16	&1.18	&\textbf{1.07}	&0.05	&0.02	&0.00	&0.99	&\textbf{0.63}	&0.73\\
&	&3	&1.00	&1.00	&1.00	&\textbf{0.88}	&1.00	&1.00	&\textbf{0.88}	&0	    &0	    &0	    &0.74	&0.61	&\textbf{0.41}\\

\cline{2-16}
&\multirow{3}{*}{100}
    &1  &1.73	&1.18	&1.34	&\textbf{1.15}	&1.36	&1.38	&\textbf{1.23}	&0.07	&0.07	&0.05	&9.73	&\textbf{3.53}	&8.66\\
&	&2	&1.55	&1.13	&1.25	&\textbf{1.07}	&1.26	&1.52	&\textbf{1.09}	&0.08	&0.09	&0.01	&11.35	&\textbf{3.33}	&5.41\\
&	&3	&1.88	&\textbf{1.08}	&1.29	&1.21	&1.33	&1.49	&\textbf{1.27}	&0.14	&0.23	&0.04	&12.53	&\textbf{4.12}	&5.45\\

\cline{2-16}
&\multirow{3}{*}{250}
    &1  &2.11	&2.11	&1.94	&\textbf{1.75}	&2.11	&2.08	&\textbf{1.84}	&0	    &0.06	&0.04	&46.08	&\textbf{38.13}	&47.03\\
&	&2	&2.30	&1.99	&1.84	&\textbf{1.70}	&2.22	&2.14	&\textbf{1.71}	&0.11	&0.17	&0.03	&54.10	&\textbf{39.93}	&47.24\\
&	&3	&2.24	&1.97	&1.79	&\textbf{1.74}	&2.14	&1.97	&\textbf{1.80}	&0.10	&0.14	&0.04	&72.40	&\textbf{38.87}	&39.59\\

\cline{2-16}
&500  &1  &2.57	&2.57	&2.13	&\textbf{2.09}	&2.57	&2.47	&\textbf{2.11}	&0	&0.18	&0.03	&277.90	&255.59	&\textbf{144.00}\\

\hline

\multirow{10}{*}{15}	&\multirow{3}{*}{50}
    &1  &1.07	&1.01	&0.98	&\textbf{0.92}	&1.06	&1.06	&\textbf{0.93}	&0.02	&0.03	&0.01	&0.89	&\textbf{0.69}	&1.25\\
&	&2	&0.99	&0.99	&0.99	&\textbf{0.88}	&0.99	&0.99	&\textbf{0.89}	&0	    &0	    &0.01	&0.81	&\textbf{0.67}	&1.13\\
&	&3	&0.89	&0.84	&0.89	&\textbf{0.79}	&0.88	&0.89	&\textbf{0.81}	&0.01	&0	    &0.02	&0.92	&\textbf{0.66}	&1.51\\

\cline{2-16}
&\multirow{3}{*}{100}
    &1  &1.17	&1.04	&1.07	&\textbf{0.97}	&1.08	&1.16	&\textbf{0.97}	&0.02	&0.03	&0	    &9.81	&3.65	&\textbf{3.51}\\
&	&2	&1.14	&1.02	&0.97	&\textbf{0.93}	&1.06	&0.99	&\textbf{0.94}	&0.04	&0.05	&0.01	&8.44	&4.35	&\textbf{3.47}\\
&	&3	&1.39	&\textbf{0.91}	&1.08	&0.99	&1.06	&1.33	&\textbf{1.01}	&0.09	&0.12	&0.01	&15.93	&3.71	&\textbf{2.88}\\

\cline{2-16}
&\multirow{3}{*}{250}
    &1  &2.05	&1.33	&1.34	&1.26	&1.46	&1.46	&\textbf{1.33}	&0.07	&0.07	&0.07	&140.88	&\textbf{52.83}	&93.75\\
&	&2	&2.08	&\textbf{1.28}	&1.31	&1.41	&\textbf{1.39}	&1.64	&1.41	&0.07	&0.29	&0.00	&165.90	&\textbf{45.20}	&65.31\\
&	&3	&1.71	&1.28	&1.23	&\textbf{1.07}	&1.38	&1.62	&\textbf{1.09}	&0.06	&0.16	&0.02	&128.30	&41.15	&\textbf{22.74}\\

\cline{2-16}
&500  &1  &2.13	&\textbf{1.64}	&1.77	&1.66	&1.87	&1.89	&\textbf{1.68}	&0.08	&0.09	&0.03	&884.69	&355.67	&\textbf{271.71}\\

\hline

\multirow{3}{*}{20}	

&100  &1  &0.98	&0.94	&0.98	&\textbf{0.83}	&0.98	&0.98	&\textbf{0.84}	&0.01	&0	    &0.01	&3.98	&\textbf{3.11}	&4.06\\
&250  &1  &1.17	&1.17	&1.17	&\textbf{0.98}	&1.17	&1.17	&\textbf{1.01}	&0	    &0	    &0.02	&53.33	&39.26	&\textbf{20.77}\\
&500  &1  &2.06	&1.35	&1.26	&\textbf{1.11}	&1.53	&1.86	&\textbf{1.14}	&0.07	&0.32	&0.02	&843.71	&314.09	&\textbf{268.58}\\

\hline

\multirow{3}{*}{25}	
&100  &1  &0.88	&0.88	&0.84	&\textbf{0.80}	&0.88	&0.88	&\textbf{0.82}	&0	    &0.01	&0.02	&4.70	&3.77	&\textbf{2.93}\\
&250  &1  &0.99	&0.99	&0.97	&\textbf{0.91}	&0.99	&0.98	&\textbf{0.91}	&0	    &0.00	&0.00	&58.07	&45.87	&\textbf{20.07}\\
&500  &1  &1.77	&1.17	&1.13	&\textbf{1.01}	&1.32	&1.63	&\textbf{1.03}	&0.05	&0.23	&0.01	&957.09	&336.30	&\textbf{220.94}\\

\hline

\end{tabular}
\medskip
\caption{Comparison of the experiment's results obtained by different heuristics.}\label{tab:res}
\end{table}

It is seen in the table that in more than in a half of all cases VNS works faster than other algorithms. But note that often, especially in small size cases, the difference in running time is not so significant. Besides, in the overwhelming majority of the cases, VNS constructs the best solution among all algorithms. Therefore, in general, the superiority of VNS is obvious. In some cases, especially when $D$ is not too large, ACO yields the better solution than VNS and GLS. But often, even when the best solution found by ACO outperforms the best solution found by VNS, the average objective value of ACO remains inferior in quality than that of VNS: for example, see the cases ($D = 7, n = 100, $ nr = $1,2$) and ($D = 10, n = 100, $ nr = $3$). Although GLS never appeared to be the best among all algorithms, it is worth to say that it often significantly improves $T_{CH}$ and builds solutions that are very close to those constructed by ACO and VNS, in terms of the objective function. Moreover, in some cases GLS outperforms one of the other algorithms: for example, see the values of $W_{best}$ in the cases ($D = 10, n = 250, $ nr = $2,3$), ($D = 15, n = 250, $ nr = $2$), and ($D = 20, n = 500, $ nr = $3$), and see the values of $W_{av}$ in the cases ($D = 7, n = 100, $ nr = $3$), ($D = 7, n = 250, $ nr = $2,3$), and ($D = 10, n = 250, $ nr = $2,3$).

In some cases, both algorithms ACO and GLS failed to improve the initial solution $T_{CH}$. This fact can be explained in the following way. These heuristics don't apply local search procedure directly to the initial solution, but they improve the derivative results of the initial solutions: mutated offspring or decoded ant path. Most probably, when solution space is rather large, there exists a risk that these two algorithms will explore only the solutions that are worse than initial solution, which may not lead to its improvement. It would be helpful to look into these cases deeper and to analyse the behaviour of the algorithms on them. Anyway, we believe that both algorithms ACO and GLS have a potential to be improved. It is also worth to say that in some cases the initial solution was significantly improved. For example, see the case ($D = 20, n = 500, $ nr = $1$), where the initial solution was improved almost twice. In particular, this gives us the following negative result regarding the constructive heuristics from \cite{PloEr19}: none of them provides an approximation with a guaranteed factor less than $2.06/1.11 \approx 1.856$.

\begin{figure}[!hbtp]
\centering
\subfloat[\label{fig:expres_CH}MPIR. $W(T) = 2.06$]{\includegraphics[width=0.5\textwidth]{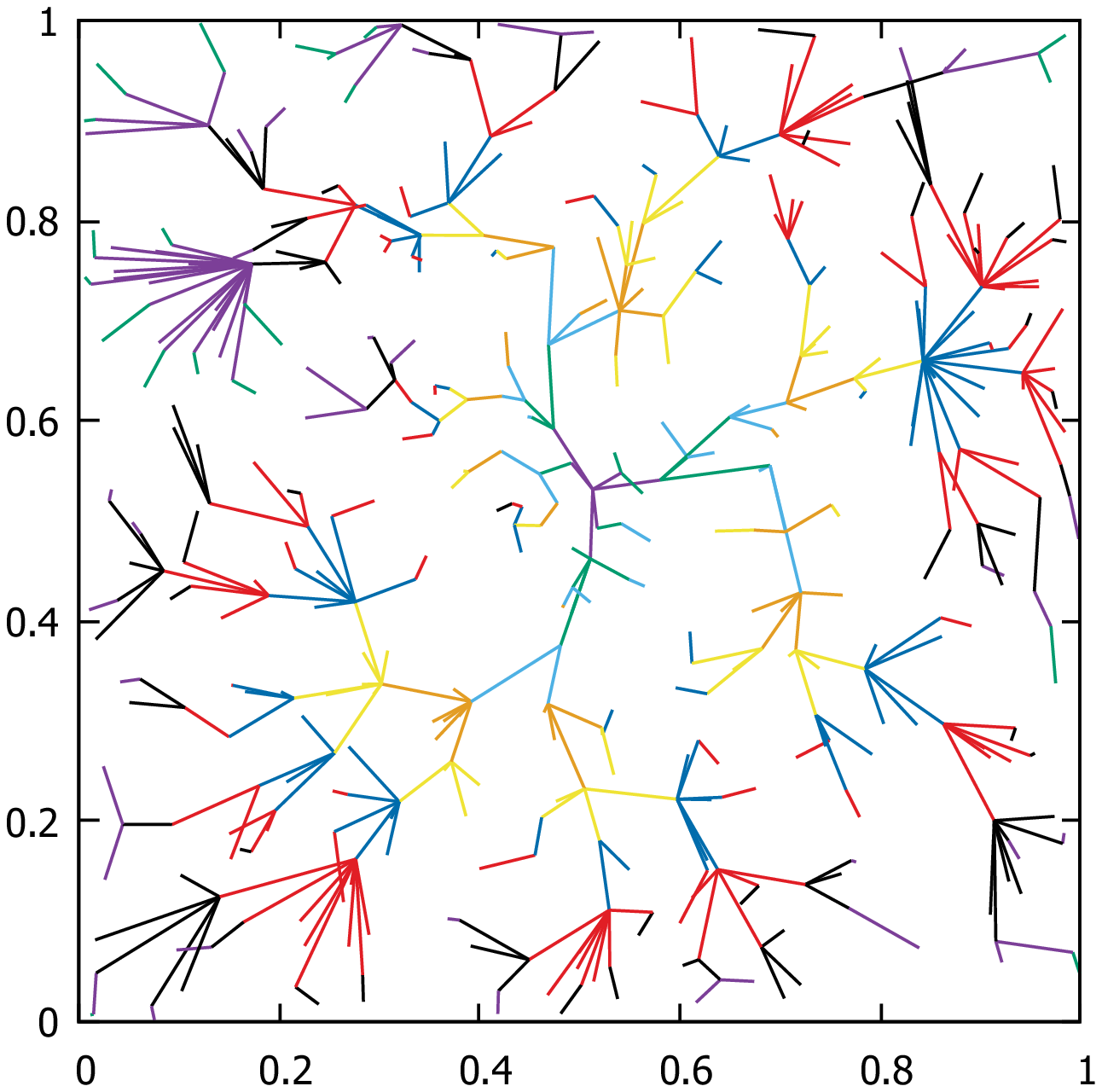}} \hfill
\subfloat[\label{fig:expres_VNS}VNS. $W(T) = 1.11$]{\includegraphics[width=0.5\textwidth]{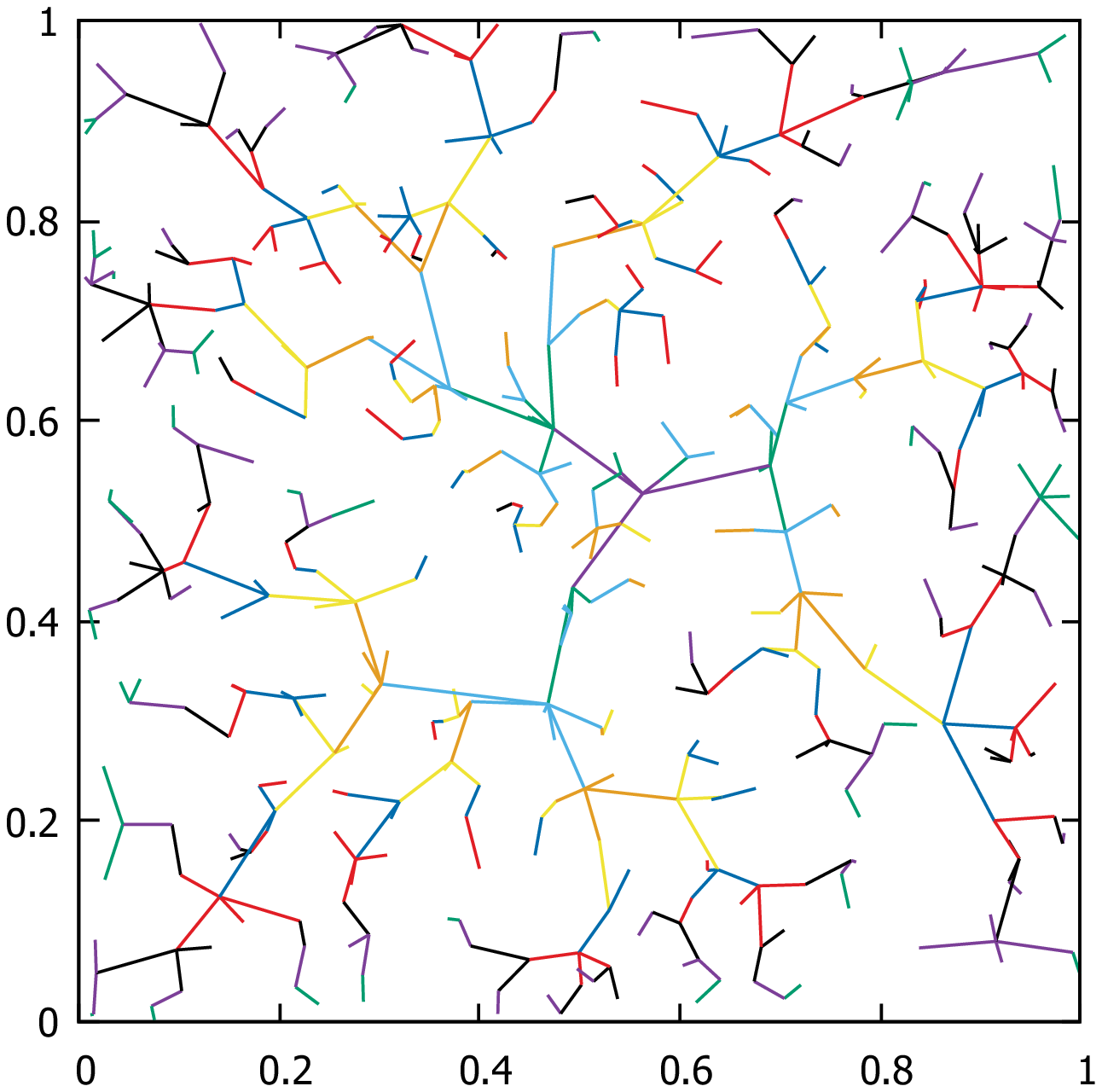}} \hfill
\subfloat[\label{fig:expres_GLS}GLS. $W(T) = 1.26$]{\includegraphics[width=0.5\textwidth]{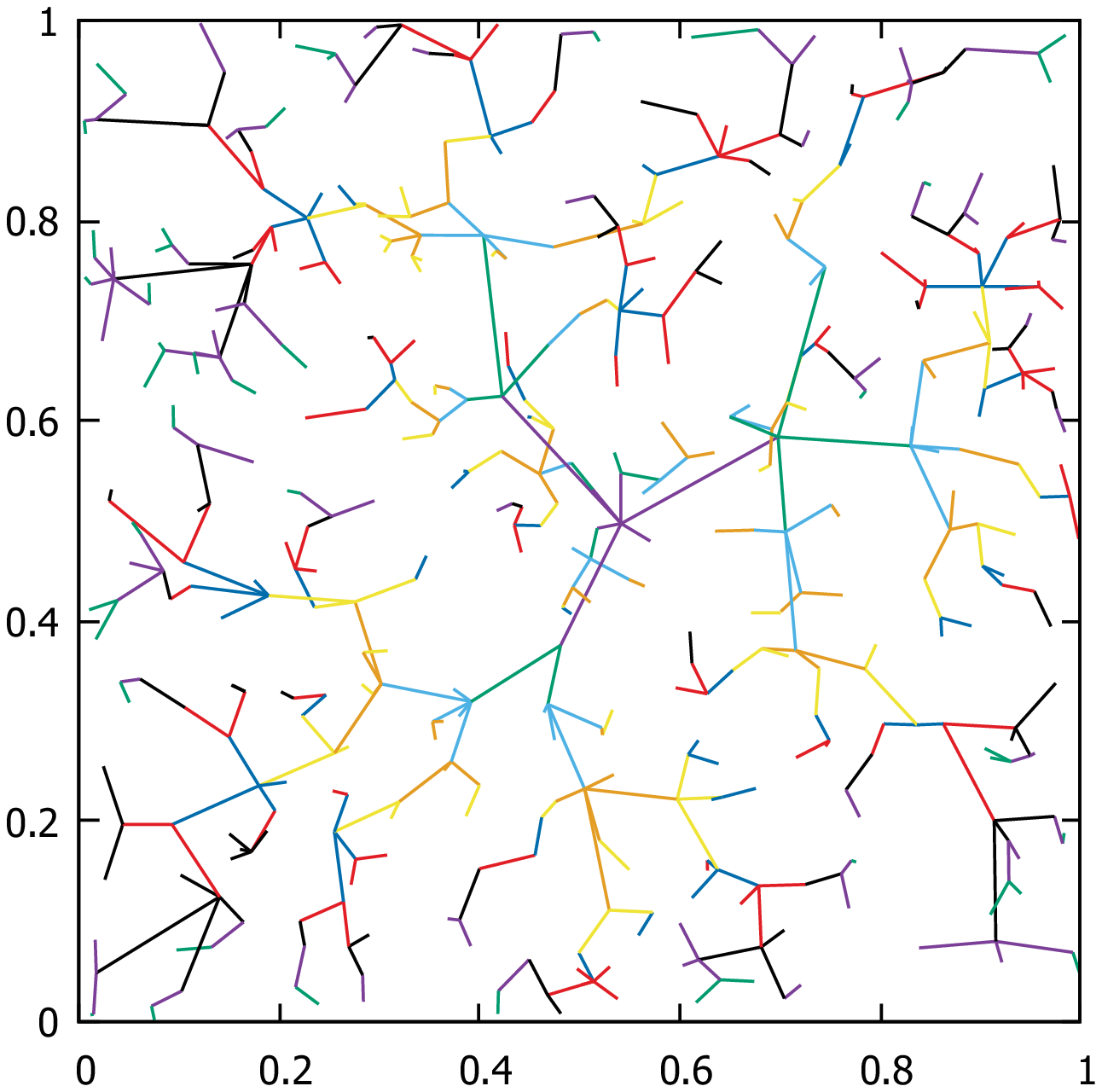}} \hfill
\subfloat[\label{fig:expres_ACO}ACO. $W(T) = 1.35$]{\includegraphics[width=0.5\textwidth]{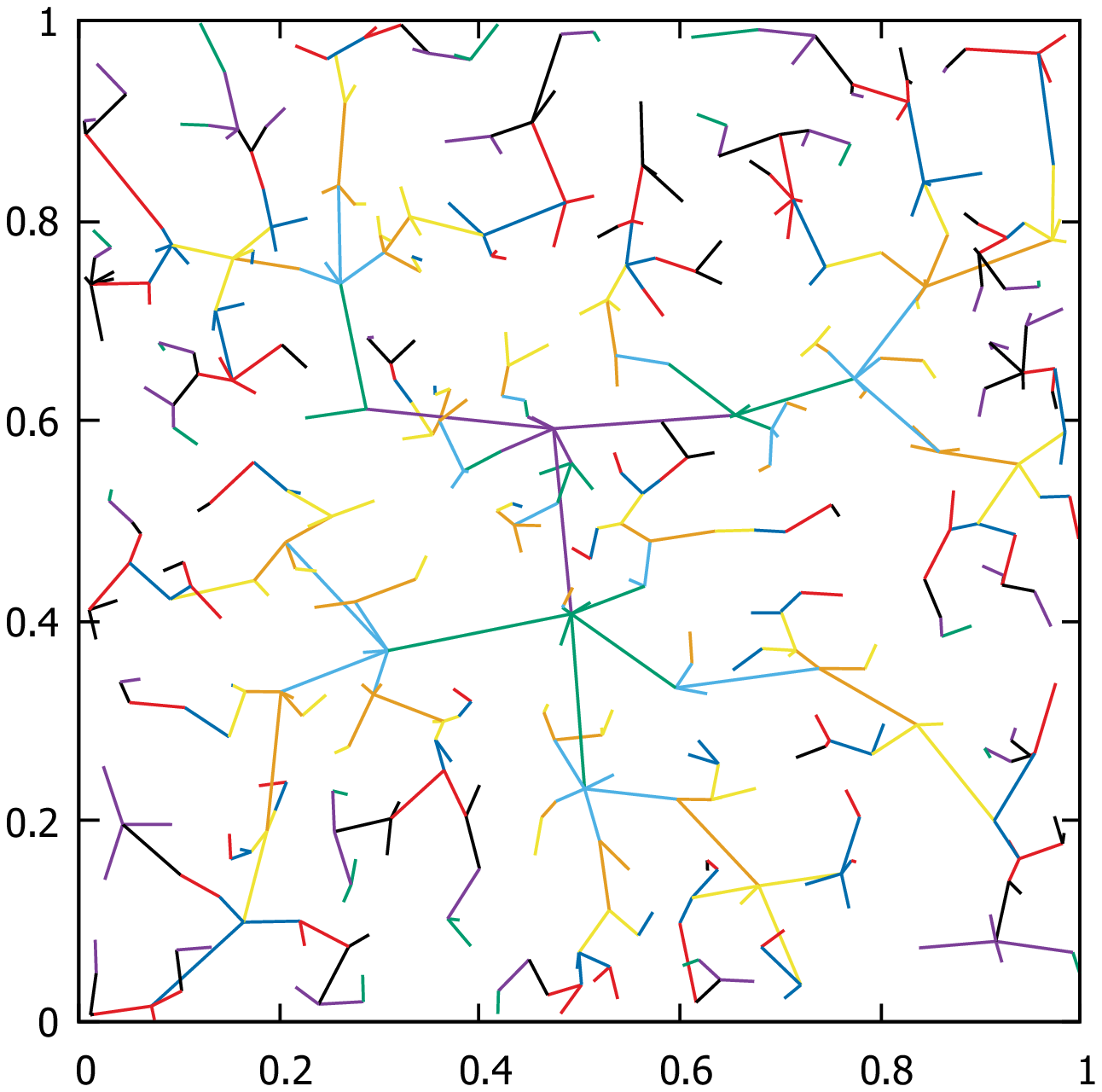}} \hfill
\caption{Best algorithms results on the same instance. $D = 20, n = 500$} \label{fig:expres}
\end{figure}

As an illustration, we also present in Fig. \ref{fig:expres} the best solutions that were obtained by different algorithms on the same instance when $D = 20, n = 500$. We chose this case, because of the big gap between the constructive heuristics and the metaheuristics results, that was discussed in the previous paragraph. For the convenience, the edges that remote from a center by an equal distance (i.e., hops count) are colored in the same color. This helps to easily verify that each tree is feasible, since the hops bound is never violated. Since the diameter bound is even in this case, there is the only center in all trees. In this case MPIR constructed the best solution among all constructive heuristics from \cite{PloEr19}. The difference between the tree constructed by MPIR, and the new metaheuristic algorithms, is seen. In the solution obtained by the constructive heuristic MPIR, a part of a tree that lies far away from the center has a star-like structure, which is not desirable, because in this case a lot of rather long edges are connected with "star centers" (i.e., the vertices with high degree). Note that the trees constructed by VNS, GLS, and ACO, have no such star-like parts: the longer edges at the backbone allow to get rid of the need for the vertices with high degree.

\section{Conclusion} \label{sC}
In this paper, we considered the NP-hard Min-Power Bounded-Hops Symmetric Connectivity Problem. For its approximation solution, we proposed three different heuristic algorithms that are based on such known metaheuristics as variable neighborhood search, ant colony optimization, and genetic local search. To the best of our knowledge, this is the first application of such kind of heuristics to this problem. We implemented all the proposed algorithms and conducted the numerical experiment on different test instances that were generated on the data sets taken from the Beasley's OR-Library. The simulation showed that, in general, our methods allow to significantly improve the solution built by the best known polynomial constructive heuristics. In most cases, VNS based heuristic appeared to be more efficient than other methods both in terms of objective function and running time.


\begin{thebibliography}{20}

\bibitem{R96}
Rappaport, T.S.:
Wireless Communications: Principles and Practices. Prentice Hall (1996)

\bibitem{CPS99}
Clementi, A.E.F., Penna, P., Silvestri, R.:
Hardness Results for the Power Range
Assignment Problem in Packet Radio Networks. In: Hochbaum, D.S., Jansen, K.,
Rolim, J.D.P., Sinclair, A. (eds.) RANDOM 1999 and APPROX 1999. LNCS,
vol. 1671, Springer, Heidelberg 197–-208 (1999)

\bibitem{KKKP00}
Kirousis, L.M., Kranakis, E., Krizanc, D., Pelc, A.:
Power Consumption in Packet
Radio Networks. Theoretical Computer Science 243(1--2), 289-–305 (2000)

\bibitem{CMZ02}
Calinescu, G., Mandoiu, I.I., Zelikovsky, A.:
Symmetric Connectivity with Minimum
Power Consumption in Radio Networks. In: Baeza-Yates, R.A., Montanari,
U., Santoro, N. (eds.) Proc. 2nd IFIP International Conference on Theoretical
Computer Science. IFIP Conference Proceedings, vol. 223, Kluwer,
Dordrecht 119–-130 (2002)

\bibitem{CNSCL03}
Cheng, X., Narahari, B., Simha, R., Cheng, M.X., Liu, D.:
Strong Minimum Energy Topology in Wireless Sensor Networks: NP-Completeness and Heuristics. IEEE
Transactions on Mobile Computing 2(3), 248–-256 (2003)

\bibitem{CN02}
Chu, T., Nikolaidis, I.:
Energy efficient broadcast in mobile ad hoc networks. In Proc. AD-HOC Networks and
Wireless (2002)

\bibitem{PS06}
Park, J., Sahni, S.:
Power Assignment For Symmetric Communication In Wireless
Networks. In: Proceedings of the 11th IEEE Symposium on Computers and
Communications (ISCC), Washington, pp. 591–596. IEEE Computer Society, Los
Alamitos (2006)

\bibitem{ACMPTZ06}
Althaus, E., Calinescu, G., Mandoiu, I.I., Prasad, S.K., Tchervenski, N., Zelikovsky,
A.: Power Efficient Range Assignment for Symmetric Connectivity in Static Ad Hoc
Wireless Networks. Wireless Networks 12(3), 287-–299 (2006)

\bibitem{EPS13}
Erzin, A., Plotnikov, R., Shamardin, Y.: On some polynomially solvable cases
  and approximate algorithms in the optimal communication tree construction
  problem, Journal of Applied and Industrial Mathematics 7 142--152 (2013)

\bibitem{WM09}
Wolf, S., Merz, P.:
Iterated Local Search for Minimum Power Symmetric Connectivity in Wireless Networks.
Volume 5482 of the series Lecture Notes in Computer Science 192--203 (2009)

\bibitem{EP15}
Erzin, A., Plotnikov, R.:
Using {VNS} for the optimal synthesis of the communication tree in wireless sensor networks, Electronic Notes in Discrete Mathematics 47 21--28 (2015)


\bibitem{EMP16_COR}
Erzin, A., Mladenovic, N., Plotnikov, R.:
Variable neighborhood search variants for Min-power symmetric connectivity problem. Computers \& Operations Research 78 557 -- 563 (2017)

\bibitem{PEM18_VNDS}
Plotnikov, R., Erzin, A., Mladenovic, N. VNDS for the Min-Power Symmetric Connectivity // Optimization Letters (2018) DOI 10.1007/s11590-018-1324-0

\bibitem{Clementi00_1}
Clementi, A.E.F., Ferreira, A., Penna, P., Perennes, S., Silvestri, R.: The minimum
range assignment problem on linear radio networks. In: Paterson, M. (ed.) ESA
2000. LNCS, vol. 1879, pp. 143–154. Springer, Heidelberg (2000)

\bibitem{Calinescu06}
Calinescu, G., Kapoor, S., Sarwat, M.: Bounded-hops power assignment in ad hoc
wireless networks. Discrete Applied Mathematics 154(9), 1358–-1371 (2006)

\bibitem{Carmi15}
Paz Carmi, Lilach Chaitman-Yerushalmi, and Ohad Trabelsi: On the Bounded-Hop Range
Assignment Problem. F. Dehne et al. (Eds.): WADS 2015, LNCS 9214, pp. 140-–151 (2015)

\bibitem{Clementi00}
A.E.F. Clementi, P. Penna, and R. Silvestri. On the power assignment problem in
radio networks. Electronic Colloquium on Computational Complexity (ECCC),
(054), 2000

\bibitem{PloEr19}
R. Plotnikov and A. Erzin. Constructive Heuristics for Min-Power Bounded-Hops Symmetric Connectivity Problem (2019) https://arxiv.org/abs/1902.06796

\bibitem{Gruber06}
M.Gruber, J. van Hemert, G. R. Raidl. Neighbourhood searches for the bounded diameter minimum spanning tree problem embedded in a VNS, EA, and ACO. In Proc. 8th annual conference on Genetic and evolutionary computation (GECCO'06) Seattle, Washington, USA, July 08--12, (2006), pp. 1187--1194. DOI:10.1145/1143997.1144185

\bibitem{HansenVNS}
Hansen, P., Mladenovic, N.:
Variable neighborhood search: Principles and applications, European Journal of Operational Research 130 449--467 (2001)



\end{thebibliography}
\end{document}